\documentclass[11 pt]{article}
\title{Small generating sets for the Torelli group\vspace{-6pt}}
\author{Andrew Putman\footnote{Supported in part by NSF grant DMS-1005318}\vspace{-6pt}}
\usepackage[compact]{titlesec}
\usepackage{setspace}
\usepackage{amsmath}
\usepackage{amssymb}
\usepackage{amsthm}
\usepackage{epsfig}
\usepackage[vscale=0.723]{geometry}
\usepackage{amsfonts}
\usepackage{calc}
\usepackage{amscd}
\usepackage[font=small]{caption}
\usepackage{pinlabel}
\usepackage[all,cmtip]{xy}

\usepackage{mathpazo}  %roman=Palatino, math=Palatino where possible
\usepackage{mathptmx}  %roman=Times, math=Times where possible

\theoremstyle{plain}
\newtheorem{theorem}{Theorem}[section]
\newtheorem{maintheorem}{Theorem}

\newtheorem{lemma}[theorem]{Lemma}

\newtheorem{question}[theorem]{Question}

\newtheorem{case}{Case}

\newcommand\BeginCases{\setcounter{case}{0}}

\theoremstyle{definition}

\newtheorem*{definition}{Definition}

\theoremstyle{remark}
\newtheorem*{remark}{Remark}

\DeclareMathOperator{\Mod}{Mod}
\newcommand\Torelli{\ensuremath{{\mathcal I}}}

\DeclareMathOperator{\SL}{SL}

\newcommand\Curves{\ensuremath{\mathcal{C}}}

\newcommand\Z{\ensuremath{\mathbb{Z}}}

\DeclareMathOperator{\HH}{H}

\DeclareMathOperator{\Aut}{Aut}

\newcommand\Span[1]{\ensuremath{\langle #1 \rangle}}
\newcommand\CaptionSpace{\hspace{0.2in}}

\newcommand\Set[2]{\ensuremath{\{\text{#1 $|$ #2}\}}}

\newcommand\Figure[3]{
\begin{figure}[t]
\centering
\centerline{\psfig{file=#2,scale=60}}
\caption{#3}
\label{#1}
\end{figure}}

\newcommand\HandleGraph{\ensuremath{\mathcal{H}_{a,b}}}

\begin{document}

\maketitle

\vspace{-23pt}
\begin{abstract}
Proving a conjecture of Dennis Johnson, we show that the Torelli subgroup $\Torelli_{g}$ of the genus
$g$ mapping 
class group has a finite generating set whose size grows cubically with respect to $g$.  Our main tool
is a new space called the handle graph on which $\Torelli_g$ acts cocompactly.
\end{abstract}

\section{Introduction}
\label{section:introduction}

Let $\Sigma_{g,n}$ be a compact connected oriented genus $g$ surface with $n$ boundary components.  The
{\em mapping class group} of $\Sigma_{g,n}$, denoted $\Mod_{g,n}$, is the group of orientation-preserving
homeomorphisms of $\Sigma_{g,n}$ that fix the boundary pointwise modulo isotopies that
fix the boundary pointwise.  We will often omit the $n$ if it vanishes.  For $n \leq 1$, the {\em Torelli group}, denoted
$\Torelli_{g,n}$, is the kernel of the action of $\Mod_{g,n}$ on $\HH_1(\Sigma_{g,n};\Z)$.  The Torelli
group has been the object of intensive study ever since the seminal work of Dennis Johnson in
the early '80's.  See \cite{JohnsonSurvey} for a survey of Johnson's work.

\paragraph{Finite generation of Torelli.}
One of Johnson's most celebrated theorems says that $\Torelli_{g,n}$ is finitely generated for
$g \geq 3$ and $n \leq 1$ (see \cite{JohnsonFinite}).  This is a surprising result -- 
though $\Mod_{g,n}$ is finitely
presentable, $\Torelli_{g,n}$ is an infinite-index normal subgroup of $\Mod_{g,n}$, so 
there is no reason to hope that $\Torelli_{g,n}$ has any finiteness properties.  Moreover,
McCullough and Miller \cite{McCulloughMiller} proved that $\Torelli_{2,n}$ is {\em not} finitely generated for $n \leq 1$,
and later Mess \cite{MessThesis} proved that $\Torelli_{2}$ is an infinite rank free
group.

\paragraph{Johnson's generating set.}
Johnson's generating set for $\Torelli_{g,n}$ when $g \geq 3$ and $n \leq 1$ is enormous.  Indeed,
for $\Torelli_g$ (resp.\ $\Torelli_{g,1}$), it contains
$9 \cdot 2^{2g-3}-4g^2+2g-6$ (resp.\ $9 \cdot 2^{2g-3}-4g^2+4g-5$) 
elements.  In \cite{JohnsonAbelianization}, Johnson
proved that the abelianization of $\Torelli_g$ (resp.\ $\Torelli_{g,1}$) has rank
$\frac{1}{3}(4g^3+5g+3)$ (resp.\ $\frac{1}{3}(4g^3-g)$).  These give large lower bounds
on the size of generating sets for $\Torelli_{g,n}$; however, there is a huge gap between
this cubic lower bound and Johnson's exponentially growing generating set.  At the end
of \cite{JohnsonFinite} and in \cite[p.\ 168]{JohnsonSurvey}, Johnson conjectures 
that there should be a generating set
for $\Torelli_{g,n}$ whose size grows cubically with respect to the genus.  
Later, in \cite[Problem 5.7]{FarbProblems} Farb asked whether there at least exists a generating set
whose size grows polynomially.

\paragraph{Main theorem.}
In this paper, we prove Johnson's conjecture.  Our main theorem is as follows.

\begin{maintheorem}
\label{theorem:maincount}
For $g \geq 3$, the group $\Torelli_g$ has a generating set of size at most $57 \binom{g}{3}$ and the
group $\Torelli_{g,1}$ has a generating set of size at most $57 \binom{g}{3} + 2g + 1$.
\end{maintheorem}

\Figure{figure:main}{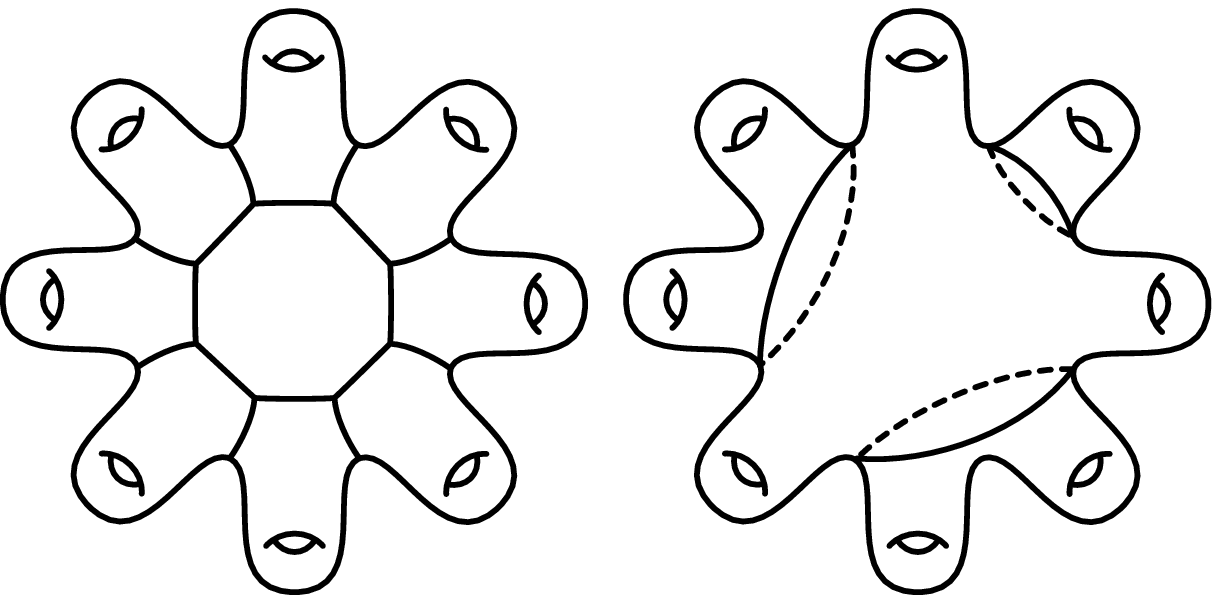}{a. The subsurfaces $R_i' \cong \Sigma_{1,1}$.  To avoid
cluttering the picture, the portion of the boundaries of the $R_i'$ which lie on the back side
the figure are not drawn. \CaptionSpace b. A subsurface isotopic to $R_{136}$.}

\noindent
The generating set we construct was conjectured to generate $\Torelli_{g,n}$ by Brendle
and Farb \cite{BrendleFarb}.  To describe it,
we must introduce some notation.  As in Figure \ref{figure:main}.a, let $R_1',\ldots,R_g'$
be $g$ subsurfaces of $\Sigma_g$ each homeomorphic to $\Sigma_{1,1}$ such that the following
hold.  Interpret all indices modulo $g$.
\begin{itemize}
\item If $1 \leq i < j \leq g$ satisfy $i \notin \{j-1,j+1\}$, then $R_i' \cap R_j' = \emptyset$.
\item For all $1 \leq i \leq g$, the intersection $R_i' \cap R_{i+1}'$ is homeomorphic to an
interval.  
\end{itemize}
For $1 \leq i < j < k \leq g$, define a subsurface $R_{ijk}$ of $\Sigma_g$ by
$R_{ijk} = \overline{\Sigma_g \setminus \bigcup_{l \neq i,j,k} R_l'}$.
Thus $R_{ijk}$ is a genus $3$ surface with at most $3$ boundary components such that
$R_i',R_j',R_k' \subset R_{i,j,k}$ (see Figure \ref{figure:main}.b).

If $S$ is a subsurface of $\Sigma_g$, define $\Mod(\Sigma_g,S)$ to be the subgroup of $\Mod_g$
consisting of mapping classes that can be realized by homeomorphisms supported on $S$
and $\Torelli(\Sigma_g,S)$ to equal $\Torelli_g \cap \Mod(\Sigma_g,S)$.  The key result for the 
proof of Theorem \ref{theorem:maincount} is the following theorem.

\begin{maintheorem}
\label{theorem:main}
For $g \geq 3$, the group $\Torelli_g$ is generated by the set
$\bigcup_{1 \leq i < j < k \leq g} \Torelli(\Sigma_g,R_{ijk})$.
\end{maintheorem}

\noindent
Using Johnson's work, it is easy to see that $\Torelli(\Sigma_g,R_{ijk})$ is finitely generated by a generating set with at most $57$
generators (see Lemma \ref{lemma:rijkgen}).  Also, standard techniques (see Lemma \ref{lemma:addboundarycpt}) 
show that if $\Torelli_g$ has a generating set with $k$ elements,
then $\Torelli_{g,1}$ has a generating set with $k + 2g + 1$ elements.  Since there are $\binom{g}{3}$ subsurfaces
$R_{ijk}$, Theorem \ref{theorem:maincount} follows from Theorem \ref{theorem:main}.

\begin{remark}
To illustrate the relative sizes of our generating sets, Johnson's generating
set for $\Torelli_{20}$ contains more than one trillion elements while our generating
set for $\Torelli_{20}$ has $64980$ elements.
\end{remark}

\paragraph{New proof of Johnson's theorem.}
Our deduction of Theorem \ref{theorem:maincount} from Theorem \ref{theorem:main} depends
on Johnson's theorem that $\Torelli_3$ is finitely generated.
However, Hain \cite{HainTorelli3} 
has recently announced a direct conceptual proof that $\Torelli_3$ is finitely
generated.
Hain's proof uses special properties of the moduli space of genus $3$ Riemann surfaces and
cannot be easily generalized to $g > 3$.  Combining this with our paper, we obtain a
new proof that $\Torelli_{g,n}$ is finitely generated
for $g \geq 3$ and $n \leq 1$.

Our new proof is more conceptual than Johnson's original one.  To illustrate this, we
will sketch Johnson's proof.  He starts by writing down an enormous finite subset $S \subset \Torelli_{g,n}$
which is known (from work of Powell \cite{PowellTorelli}) to normally generate $\Torelli_{g,n}$ as
a subgroup of $\Mod_{g,n}$.  Letting $T$ be a standard generating set for $\Mod_{g,n}$, Johnson
then proves via a laborious computation that for $t \in T$ and $s \in S$, the element
$t s t^{-1} \in \Torelli_{g,n}$ can be written as a word in $S$.  This implies
that the subgroup $\Gamma$ of $\Torelli_{g,n}$ generated by $S$ is a normal subgroup of $\Mod_{g,n}$,
and thus that $\Gamma = \Torelli_{g,n}$.

\begin{remark}
Our proof of Theorem \ref{theorem:main} appeals to a theorem of \cite{PutmanConnectivityNote} whose
proof depends on Johnson's theorem.  However, Hatcher and Margalit \cite{HatcherMargalit} have
recently given a new proof of this result that is independent of Johnson's work.
\end{remark}

\paragraph{Nature of generators.}
Some basic elements of $\Torelli_{g,n}$ are as follows (see, e.g., \cite{PutmanCutPaste}).  
If $x$ is a simple closed curve on $\Sigma_{g,n}$,
then denote by $T_x \in \Mod_{g,n}$ the Dehn twist about $x$.  If $x$ is a separating simple closed curve, then
$T_x \in \Torelli_{g,n}$; these are called {\em separating twists}.  If $x$ and $y$ are disjoint
homologous nonseparating simple closed curves, then $T_x T_y^{-1} \in \Torelli_{g,n}$; these are called
{\em bounding pair maps}.  Following work of Birman \cite{BirmanSiegel}, Powell \cite{PowellTorelli}
proved that $\Torelli_{g,n}$ is generated by bounding pair maps and separating twists for $g \geq 1$
and $n \leq 1$ (see \cite{PutmanCutPaste} and \cite{HatcherMargalit} for alternate proofs).  
Johnson's finite generating set for $\Torelli_{g,n}$ for $g \geq 3$ and $n \leq 1$
consists entirely of bounding pair maps.  It follows easily from our proofs
of Lemma \ref{lemma:addboundarycpt} and \ref{lemma:rijkgen} that our generating set consists
of bounding pair maps and separating twists; see the remark after Lemma \ref{lemma:rijkgen}.

\paragraph{The handle graph.}
Our proof of Theorem \ref{theorem:main} is topological.  To prove that a group $G$ is finitely generated,
it is enough to find a connected simplicial complex upon which $G$ acts cocompactly with finitely
generated stabilizers.  We use a variant on the curve complex.  If $\gamma$ is an oriented
simple closed curve on $\Sigma_g$, then denote by $[\gamma] \in \HH_1(\Sigma_g;\Z)$ its homology class.
Also, if $\gamma_1$ and $\gamma_2$ are isotopy classes of simple closed curves on $\Sigma_g$, then
denote by $i_g(\gamma_1,\gamma_2)$ their {\em geometric intersection number}, i.e.\ the minimal possible
number of intersections between two curves in the isotopy classes of $\gamma_1$ and $\gamma_2$.  Finally,
denote by $i_a(\cdot,\cdot)$ the algebraic intersection pairing on $\HH_1(\Sigma_g;\Z)$.

\begin{definition}
Let $a,b \in \HH_1(\Sigma_g;\Z)$ satisfy $i_a(a,b)=1$.  The {\em handle graph} associated to $a$ and $b$,
denoted $\HandleGraph$, is the graph whose vertices are isotopy classes of oriented
simple closed curves on $\Sigma_g$ that are homologous to either $a$ or $b$ and where two
vertices $\gamma_1$ and $\gamma_2$ are joined by an edge exactly when $i_g(\gamma_1,\gamma_2)=1$.
\end{definition}

\noindent
We will show that $\HandleGraph / \Torelli_g$ consists of a single edge (see Lemma \ref{lemma:handlegraphquotient})
and that $\HandleGraph$ is connected for $g \geq 3$ (see Lemma \ref{lemma:handlegraphconnected}).

\paragraph{A complication.}
It would appear that we have all the ingredients in place to use the space $\HandleGraph$
to prove that $\Torelli_g$ is finitely generated.
However, there is one remaining complication.  Namely, we do not know the
answer to the following question.

\begin{question}
\label{question:stabilizerfg}
For some $g \geq 4$, let $\gamma$ be the isotopy class of a nonseparating simple closed curve on $\Sigma_g$.
Is the stabilizer subgroup $(\Torelli_g)_{\gamma}$ of $\gamma$ finitely generated?
\end{question}

\noindent
In other words, we do not know if the vertex stabilizer subgroups of the action
of $\Torelli_g$ on $\HandleGraph$ are finitely generated.  Nonetheless, in
\S \ref{section:birmangenerators} we will prove a weaker statement that suffices
to prove Theorem \ref{theorem:main}.  The proof of Theorem \ref{theorem:main} is in \S \ref{section:main}.

\paragraph{Smaller generating sets.}
A positive answer to Question \ref{question:stabilizerfg} would likely lead to a smaller
generating set for $\Torelli_g$, though of course this depends on the nature of
the finite generating sets for the stabilizer subgroups.  Let us describe one
way this could work.  For $g \geq 3$, let $\sigma_g$ be the smallest cardinality of a
generating set for $\Torelli_g$.  Consider $g \geq 4$, and fix an edge $\{\alpha,\beta\}$ of $\HandleGraph$.  The proof
of Theorem \ref{theorem:main} shows that $\Torelli_g$ is generated by $(\Torelli_g)_{\alpha} \cup (\Torelli_g)_{\beta}$.
Let $S$ be a subsurface of $\Sigma_g$ such that $S \cong \Sigma_{g-1,1}$ and $\alpha \cup \beta \subset \Sigma_g \setminus S$.
We have $\Torelli(\Sigma_g,S) \cong \Torelli_{g-1,1}$ (see \S \ref{section:torellibdry}) and
$\Torelli(\Sigma_g,S) \subset (\Torelli_g)_{\alpha}$ and $\Torelli(\Sigma_g,S) \subset (\Torelli_g)_{\beta}$.  Assume that there exists
a finite set $V_{\alpha}$ (resp.\ $V_{\beta}$) such that $(\Torelli_g)_{\alpha}$ (resp.\ $(\Torelli_g)_{\beta}$)
is generated by $\Torelli(\Sigma_g,S) \cup V_{\alpha}$ (resp.\ $\Torelli(\Sigma_g,S) \cup V_{\beta}$).  The
group $\Torelli_g$ is then generated by $\Torelli(\Sigma_g,S) \cup V_{\alpha} \cup V_{\beta}$.  Lemma
\ref{lemma:addboundarycpt} says that $\Torelli(\Sigma_g,S) \cong \Torelli_{g-1,1}$ can be generated
by $\sigma_{g-1} + 2g + 1$ elements.  Moreover, it seems likely that there exists some relatively small $K$
such that $|V_{\alpha}|,|V_{\beta}| \leq K g^2$.  This would imply that
$$\sigma_g \leq \sigma_{g-1} + 2g+1+2Kg^2.$$
Iterating this, we would get that
$$\sigma_g \leq \sigma_3 + \sum_{i=4}^g (2i+1+2Ki^2)$$
for $g \geq 4$.  This bound is cubic in $g$ (as it needs to be), but as long as $K$ is not too large
it is much smaller than $57 \binom{g}{3}$.

\paragraph{Finite presentability.}
Perhaps the most important open question about the combinatorial group theory of $\Torelli_g$ is whether or
not it is finitely presentable for $g \geq 3$.  One way of proving that a group $G$ is finitely presentable
is to construct a simply-connected simplicial complex $X$ upon which $G$ acts cocompactly with finitely presentable
stabilizer subgroups (see, e.g.,\ \cite{BrownPresentation}).  For example, Hatcher and Thurston use this
technique in \cite{HatcherThurston} to prove that the mapping class group is finitely presentable.

The handle graph $\HandleGraph$ appears to be the first example of a useful space upon which $\Torelli_g$ acts 
cocompactly (of course, there are trivial non-useful examples of such spaces; for example, 
the Cayley graph of $\Torelli_g$ or a $1$-point space).  Unfortunately, while $\HandleGraph$ is connected
for $g \geq 3$, it is not simply connected.  Indeed, it does not even have any $2$-cells (and is not a tree).  
However, one could probably attach $2$-cells to $\HandleGraph$ to obtain a simply connected complex upon
which $\Torelli_g$ acts cocompactly.  This would not be enough, however -- one would also have to 
prove that the simplex stabilizer subgroups were finitely presentable.  In other words, this complex
would provide the inductive step in a proof that $\Torelli_g$ was finitely presentable, but one would
still need a base case. 

\paragraph{A complex that does not work.}
We close this introduction by discussing an approach to Theorem \ref{theorem:main} that does not work.
One might think of trying to prove Theorem \ref{theorem:main} using the following complex.
Let $a \in \HH_1(\Sigma_g;\Z)$ be a primitive vector.  Define $\Curves_a$ to be the graph whose
vertices are isotopy classes of oriented simple closed curves $\gamma$ on $\Sigma_g$ such that $[\gamma] = a$ and
where two vertices $\gamma$ and $\gamma'$ are joined by an edge if $i_g(\gamma,\gamma')=0$.
It is known (\cite[Theorem 1.9]{PutmanConnectivityNote}; see \cite{HatcherMargalit} for an
alternate proof) that $\Curves_a$ is connected for
$g \geq 3$.  Moreover, $\Torelli_g$ acts transitively on the vertices of $\Curves_a$.  However, it does
{\em not} act cocompactly; indeed, there are infinitely many edge orbits.  To see this, consider
edges $e_1=\{\gamma_1,\gamma_1'\}$ and $e_2=\{\gamma_2,\gamma_2'\}$ of $\Curves_a$.  Assume that there exists
some $f \in \Torelli_g$ such that $f(e_1)=e_2$.  Since $\gamma_1$ is homologous to $\gamma_1'$,
the multicurve $\gamma_1 \cup \gamma_1'$ divides $\Sigma_g$ into two subsurfaces $S_1$ and $S_1'$.  Similarly,
$\gamma_2 \cup \gamma_2'$ divides $\Sigma_g$ into two subsurfaces $S_2$ and $S_2'$.  Relabeling if necessary,
we have $f(S_1)$ isotopic to $S_2$ and $f(S_1')$ isotopic to $S_2'$.  
Since $f \in \Torelli_g$, the images of $\HH_1(S_1;\Z)$ and
$\HH_1(S_2;\Z)$ in $\HH_1(\Sigma_g;\Z)$ must be the same, and similarly for $\HH_1(S_1';\Z)$ and
$\HH_1(S_2';\Z)$.  It is easy to see that infinitely many such images occur for different
edges of $\Curves_a$, so there must be infinitely many edges orbits.  We remark
that Johnson proved in \cite[Corollary to Lemma 9 on p.\ 250]{JohnsonConjugacy} that
the images of $\HH_1(S_1;\Z)$ and $\HH_1(S_1';\Z)$ in $\HH_1(\Sigma_g;\Z)$ are a complete invariant
for the edge orbits.

\paragraph{Acknowledgments.}
I wish to thank Tara Brendle, Benson Farb, and Dan Margalit for their help.  I also wish
to thank an anonymous referee for a very helpful referee report.

\section{The Torelli group on subsurfaces}
\label{section:torellibdry}

We will need to understand how the Torelli group restricts to subsurfaces.  For a
general discussion of this, see \cite{PutmanCutPaste}.  In this section, we will
extract from \cite{PutmanCutPaste} results on two kinds of subsurfaces.  In
\S \ref{section:rijk}, we will show how to analyze subsurfaces like the subsurfaces
$R_{ijk}$ from \S \ref{section:introduction}.  In \S \ref{section:curvestab}, we
will show how to analyze stabilizers of nonseparating simple closed curves (which
are supported on the subsurface obtained by taking the complement of a regular neighborhood
of the curve).

\subsection{Analyzing the subsurfaces $R_{ijk}$}
\label{section:rijk}

We begin by defining groups $\Torelli_{g,n}$ for $n \geq 2$.
There is a map $\Mod_{g,n} \rightarrow \Mod_g$ induced by gluing discs to 
the boundary components of $\Sigma_{g,n}$ and extending homeomorphisms by the identity.
Define $\Torelli_{g,n}$ to be the kernel
of the resulting action of $\Mod_{g,n}$ on $\HH_1(\Sigma_g;\Z)$.  For the case $n=1$, the map
$\HH_1(\Sigma_{g,1};\Z) \rightarrow \HH_1(\Sigma_g;\Z)$ is an isomorphism, so this agrees
with our previous definition of $\Torelli_{g,1}$.

\begin{remark}
In \cite{PutmanCutPaste}, the different definitions of the Torelli group on a surface with boundary
are parametrized by partitions of the boundary components.  The above definition
of $\Torelli_{g,n}$ corresponds to the discrete partition $\{\{\beta_1\},\ldots,\{\beta_n\}\}$
of the set $\{\beta_1,\ldots,\beta_n\}$ of boundary components of $\Sigma_{g,n}$.
\end{remark}

In \cite[Theorem 1.2]{PutmanCutPaste}, a version of the Birman exact sequence
is proven for the Torelli group.  For $\Torelli_{g,n}$ with $g \geq 2$, it takes the form
\begin{equation}
\label{eqn:birman}
1 \longrightarrow \pi_1(U\Sigma_{g,n}) \longrightarrow \Torelli_{g,n+1} \longrightarrow \Torelli_{g,n} \longrightarrow 1.
\end{equation}
Here $U\Sigma_{g,n}$ is the unit tangent bundle of $\Sigma_{g,n}$.  The subgroup $\pi_1(U\Sigma_{g,n})$
of $\Torelli_{g,n+1}$ is often called the ``disc-pushing subgroup'' -- the mapping class associated to 
$\gamma \in \pi_1(U\Sigma_{g,n})$ ``pushes'' a fixed boundary component around $\gamma$ while allowing
it to rotate.
The following is an immediate consequence of \eqref{eqn:birman} 
and the fact that $\pi_1(U\Sigma_{g})$ can be generated
by $2g+1$ elements.

\begin{lemma}
\label{lemma:addboundarycpt}
$\Torelli_{g,1}$ can be generated by $k+2g+1$ elements if $\Torelli_g$ can be generated by $k$ elements.
\end{lemma}

Now assume that $S \cong \Sigma_{h,n}$ is an embedded subsurface of $\Sigma_g$ and that all the boundary
components of $S$ are non-nullhomotopic separating curves in $\Sigma_g$.  For example, $S$ could be
one of the surfaces $R_{ijk}$ from \S \ref{section:introduction}.  Letting $\Mod(S)$ be the 
mapping class group of $S$, the induced map $\Mod(S) \rightarrow \Mod_g$ is an injection.  This gives
a natural identification of $\Mod(S)$ with $\Mod(\Sigma_g,S)$.  The group $\Torelli(\Sigma_g,S)$ is
thus naturally a subgroup of $\Mod(S) \cong \Mod_{h,n}$, and in \cite[Theorem 1.1]{PutmanCutPaste}
it is proven that $\Torelli(\Sigma_g,S) = \Torelli_{h,n}$.  Johnson \cite{JohnsonFinite}
proved that $\Torelli_3$ can be generated by $35$ elements.  Applying \eqref{eqn:birman}
repeatedly, we see that $\Torelli_{3,1}$ can be generated by $42$ elements, $\Torelli_{3,2}$ by $49$ elements,
and $\Torelli_{3,3}$ by $57$ elements.  Since $R_{ijk} \cong \Sigma_{3,k}$ with $k \leq 3$,
we obtain the following.

\begin{lemma}
\label{lemma:rijkgen}
For all $1 \leq i < j < k \leq g$, the group $\Torelli(\Sigma_g,R_{ijk})$ can be generated by $57$ elements.
\end{lemma}

\begin{remark}
It is well-known (see, e.g.,\ \cite[\S 2.1]{PutmanCutPaste}) that the mapping classes corresponding to
the generators of $\pi_1(U\Sigma_{g,n})$ used to prove Lemmas \ref{lemma:addboundarycpt} and
\ref{lemma:rijkgen} can be chosen to be bounding pair maps and separating twists.  Additionally,
Johnson's minimal-size generating set for $\Torelli_{3}$ consists entirely of bounding pair maps,
so the generating set for $\Torelli(\Sigma_g,R_{ijk})$ in Lemma \ref{lemma:rijkgen} can be
taken to consist of bounding pair maps and separating twists.
\end{remark}

\subsection{Stabilizers of nonseparating simple closed curves}
\label{section:curvestab}

Let $\gamma$ be a nonseparating simple closed curve on $\Sigma_g$.  Define $\Sigma_{g,\gamma}$
to be the result of cutting $\Sigma_g$ along $\gamma$, so $\Sigma_{g,\gamma} \cong \Sigma_{g-1,2}$.
Letting $\Mod_{g,\gamma}$ be the mapping class group of $\Sigma_{g,\gamma}$, the natural map
$\Sigma_{g,\gamma} \rightarrow \Sigma_g$ induces a map $i : \Mod_{g,\gamma} \rightarrow \Mod_g$.
Define $\Torelli_{g,\gamma} = i^{-1}(\Torelli_g)$.  The map $i$ restricts to a surjection
$\Torelli_{g,\gamma} \rightarrow (\Torelli_g)_{\gamma}$, where $(\Torelli_g)_{\gamma}$ is the stabilizer
subgroup of $\gamma$.

\begin{remark}
In the notation of \cite{PutmanCutPaste}, the group $\Torelli_{g,\gamma}$ corresponds to the
Torelli group of $\Sigma_{g-1,2}$ with respect to the ``indiscrete partition'' $\{\{\beta,\beta'\}\}$ of
the boundary components $\beta$ and $\beta'$ of $\Sigma_{g,\gamma}$.  Also, the kernel of the map
$\Torelli_{g,\gamma} \rightarrow (\Torelli_g)_{\gamma}$ is isomorphic to $\Z$ and is generated by
$T_{\beta} T_{\beta'}^{-1}$, where $T_{\beta}$ and $T_{\beta'}$ are the Dehn twists about $\beta$
and $\beta'$, respectively.
\end{remark}

\Figure{figure:birmantorelli}{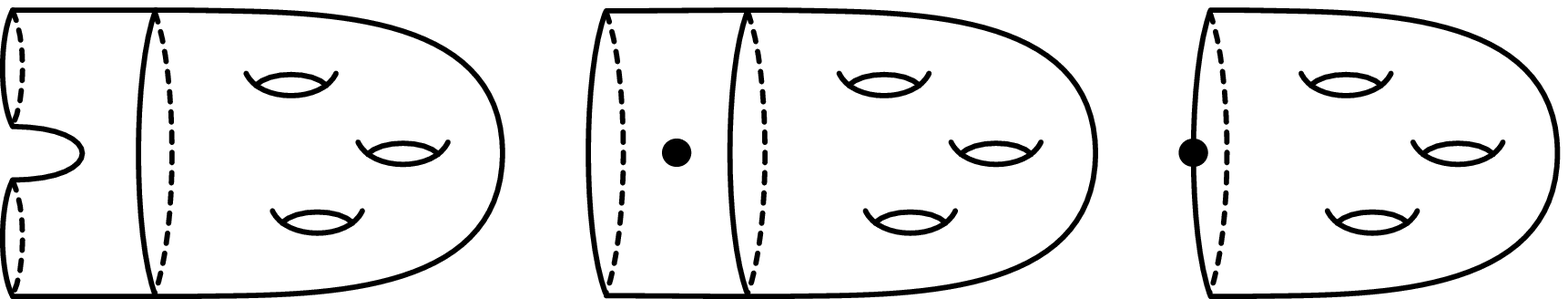}{a. The surface $\Sigma_{g,\gamma}$ and and the subsurface
$\Sigma_{g-1,1}$ of $\Sigma_{g,\gamma}$ such that the induced map 
$\Torelli_{g-1,1} \rightarrow \Torelli_{g,\gamma}$ splits the exact sequence \eqref{eqn:birmantorelli}.
\CaptionSpace b. The basepoint for $\pi_1(\Sigma_{g-1,1})$ is obtained from $\Sigma_{g,\gamma}$ by collapsing the
boundary component $\beta$ to a point. \CaptionSpace
c. The surface in b deformation retracts to $\Sigma_{g-1,1}$ such that the basepoint ends up
on the boundary component.}

In \cite[Theorem 1.2]{PutmanCutPaste}, it is proven that for $g \geq 2$ there is a short exact sequence
\begin{equation}
\label{eqn:birmantorelli}
1 \longrightarrow K_{g,\gamma} \longrightarrow \Torelli_{g,\gamma} \longrightarrow \Torelli_{g-1,1} \longrightarrow 1.
\end{equation}
Here $K_{g,\gamma} \cong [\pi_1(\Sigma_{g-1,1}),\pi_1(\Sigma_{g-1,1})]$.
This exact sequence splits via the inclusion $\Torelli_{g-1,1} \hookrightarrow \Torelli_{g,\gamma}$
induced by the inclusion $\Sigma_{g-1,1} \hookrightarrow \Sigma_{g,\gamma}$ indicated in
Figure \ref{figure:birmantorelli}.a.  In other words, the following holds.
\begin{lemma}
\label{lemma:semidirect}
$\Torelli_{g,\gamma} = K_{g,\gamma} \ltimes \Torelli_{g-1,1}$ for $g \geq 3$ and $\gamma$ a simple closed nonseparating curve on $\Sigma_g$.
\end{lemma}
\noindent
The group $\Torelli_{g-1,1}$ acts on $K_{g,\gamma} < \pi_1(\Sigma_{g-1,1})$ as follows.  As is clear
from \cite[Theorem 1.2]{PutmanCutPaste},
the basepoint for $\pi_1(\Sigma_{g-1,1})$ is as indicated in Figure \ref{figure:birmantorelli}.b.  As
shown in Figure \ref{figure:birmantorelli}.c, the surface $\Sigma_{g-1,1}$ deformation retracts
onto the surface $\Sigma_{g-1,1}$ on which $\Torelli_{g-1,1}$ is supported.  After this
deformation retract, the basepoint ends up on $\partial \Sigma_{g-1,1}$.  
Summing up, $\Torelli_{g-1,1}$ acts
on $K_{g,\gamma} < \pi_1(\Sigma_{g-1,1})$ via the action of $\Mod_{g-1,1}$ on $\pi_1(\Sigma_{g-1,1})$,
where the basepoint for $\pi_1(\Sigma_{g-1,1})$ is on $\partial \Sigma_{g-1,1}$.

\section{The handle graph is connected}
\label{section:handlegraphconnected}

In this section, we prove the following.

\begin{lemma}
\label{lemma:handlegraphconnected}
Fix $g \geq 3$.  Let $a,b \in \HH_1(\Sigma_g;\Z)$ satisfy $i_a(a,b)=1$.  Then
$\HandleGraph$ is connected.
\end{lemma}

\noindent
We will need two lemmas.  In the first, if $\epsilon$ is an oriented
arc in a surface, then $\epsilon^{-1}$ denotes the arc obtained by reversing
the orientation of $\epsilon$.

\begin{lemma}
\label{lemma:grafting}
Let the boundary components of $\Sigma_{g,2}$ be $\delta_0$ and $\delta_1$.  Choose points $v_i \in \delta_i$
for $i=0,1$ and let $\epsilon$ be an oriented properly embedded arc in $\Sigma_{g,2}$ whose
initial point is $v_0$ and whose terminal point is $v_1$.
Then for any $h \in \HH_1(\Sigma_{g,2};\Z)$, there exists an oriented properly embedded arc
$\epsilon'$ in $\Sigma_{g,2}$
whose initial point is $v_0$ and whose terminal point is $v_1$ 
such that the homology class of the loop $\epsilon' \cdot \epsilon^{-1}$
is $h$.
\end{lemma}
\begin{proof}
Gluing $(\delta_0,v_0)$ to $(\delta_1,v_1)$, we obtain a surface 
$S \cong \Sigma_{g+1}$.  Let $\alpha$ and $\ast$ be the images of $\delta_0$ and $v_0$ in $S$,
respectively.  The image
of $\epsilon$ in $S$ is an oriented simple closed curve $\beta$ with $i_g(\alpha,\beta)=1$.  There
is a natural isomorphism $\HH_1(\Sigma_{g,2};\Z) \cong [\alpha]^{\perp}$, where the orthogonal
complement is taken with respect to $i_a(\cdot,\cdot)$.  Under this identification,
we can apply \cite[Lemma A.3]{PutmanCutPaste} to find an oriented simple closed curve $\beta'$
on $S$ such that $[\beta'] = [\beta] + h$ and such that $\alpha \cap \beta' = \{\ast\}$.
Cutting $S$ open along $\alpha$, the curve $\beta'$ becomes the desired arc $\epsilon'$.
\end{proof}

\begin{lemma}
\label{lemma:connectatoa}
Let $a,b \in \HH_1(\Sigma_g;\Z)$ satisfy $i_a(a,b)=1$.  Let $\alpha_1$ and $\alpha_2$ be disjoint oriented simple
closed curves on $\Sigma_g$ such that $[\alpha_i] = a$ for $i=1,2$.  
There then exists some oriented simple closed curve $\beta$ on $\Sigma_g$ such that $[\beta] = b$ and 
$i_g(\alpha_i,\beta) = 1$ for $i=1,2$.
\end{lemma}
\begin{proof}
Let $\beta'$ be any simple closed curve on $\Sigma_g$ such that $i(\alpha_i,\beta')=1$ for $i=1,2$.  Orient
$\beta'$ so that its intersections with $\alpha_1$ and $\alpha_2$ are positive.  Let $X_1$ and $X_2$
be the two subsurfaces of $\Sigma_g$ that result from cutting $\Sigma_g$ along $\alpha_1 \cup \alpha_2$.
For $i=1,2$, the surface $X_i$ has $2$ boundary components and the intersection of $\beta'$ with $X_i$ is
an oriented properly embedded arc $\epsilon_i$ running between these boundary components.
Also, the induced map $\HH_1(X_i;\Z) \rightarrow \HH_1(\Sigma_g;\Z)$ is an injection, and we will identify
$\HH_1(X_i;\Z)$ with its image in $\HH_1(\Sigma_g;\Z)$.  The orthogonal complement
to $a$ with respect to the algebraic intersection pairing is spanned by $\HH_1(X_1;\Z) \cup \HH_1(X_2;\Z)$.
Since $i_a(a,b) = i_a(a,[\beta'])$, the homology class $b-[\beta']$ is orthogonal to $a$.  There thus
exist $h_i \in \HH_1(X_i;\Z)$ for $i=1,2$ such that $b = [\beta'] + h_1 + h_2$.  Lemma \ref{lemma:grafting}
says that for $i=1,2$ there exists an oriented properly embedded arc $\epsilon_i'$ in 
$X_i$ with the same endpoints as
$\epsilon_i$ such that the homology class of the loop 
$\epsilon_i' \cdot \epsilon_i^{-1}$ equals $h_i$.  Letting $\beta$ be
the loop $\epsilon_1' \cdot \epsilon_2'$, it follows that $[\beta] = [\beta'] + h_1 + h_2 = b$, as desired.
\end{proof}

\begin{proof}[{Proof of Lemma \ref{lemma:handlegraphconnected}}]
Let $\delta$ and $\delta'$ be
vertices of $\HandleGraph$.  We will construct a path in $\HandleGraph$ from $\delta$ to $\delta'$.  Without
loss of generality, $[\delta] = [\delta'] = a$.  By \cite[Theorem 1.9]{PutmanConnectivityNote}
(see \cite{HatcherMargalit} for an alternate proof), we can
find a sequence
$$\delta = \alpha_1,\alpha_2,\ldots,\alpha_n = \delta'$$
of isotopy classes of oriented simple closed curves on $\Sigma_g$ such that $[\alpha_i] = a$
for $1 \leq i \leq n$ and $i_g(\alpha_i,\alpha_{i+1})=0$ for $1 \leq i < n$ (this is where
we use the condition $g \geq 3$).  Lemma \ref{lemma:connectatoa} implies that there
exist isotopy classes $\beta_1,\ldots,\beta_{n-1}$ of oriented simple closed curves on $\Sigma_g$
such that $[\beta_i] = b$ and $i_g(\alpha_i,\beta_i) = i_g(\alpha_{i+1},\beta_i) = 1$ for
$1 \leq i < n$.  Since $\beta_i$ is adjacent to both $\alpha_i$ and $\alpha_{i+1}$ in $\HandleGraph$,
the desired path from $\delta$ to $\delta'$ is thus
\[\delta = \alpha_1,\beta_1,\alpha_2,\beta_2,\ldots,\beta_{n-1},\alpha_n = \delta'.\qedhere\]
\end{proof}

\section{Generating the stabilizer of a nonseparating simple closed curve}
\label{section:birmangenerators}

Let the subsurfaces $R_i'$ of $\Sigma_g$ be as in the introduction.  Define 
$S_i = \overline{\Sigma_{g} \setminus R_i'}$.  The goal of this section is to prove
the following lemma.

\begin{lemma}
\label{lemma:stabgen}
Assume that $g \geq 4$.  Let $\gamma$ be the isotopy class of a simple closed nonseparating
curve on $\Sigma_g$ that is contained in $R_1'$.  Then the subgroup $(\Torelli_{g})_{\gamma}$ of $\Torelli_g$ stabilizing $\gamma$ is contained in
the subgroup of $\Torelli_g$ generated by $\cup_{i=1}^g \Torelli(\Sigma_{g},S_i)$.
\end{lemma}

Before proving this, we need a technical lemma.
Set $\pi = \pi_1(\Sigma_{g,1},\ast)$, where $\ast \in \partial \Sigma_{g,1}$.  
Let $T_1',\ldots,T_g'$ be disjoint subsurfaces
of $\Sigma_{g,1}$ such that $T_i' \cong \Sigma_{1,1}$ and $T_i' \cap \partial \Sigma_{g,1} = \emptyset$ for $1 \leq i \leq g$ (see
Figure \ref{figure:birmangenset}.a).  Define 
$T_i = \overline{\Sigma_{g,1} \setminus T_i'}$.  We have $T_i \cong \Sigma_{g-1,2}$ and $\ast \in T_i$
for $1 \leq i \leq g$.  The maps $\pi_1(T_i,\ast) \rightarrow \pi_1(\Sigma_{g,1},\ast)$ and
$\HH_1(T_i';\Z) \rightarrow \HH_1(\Sigma_{g,1};\Z)$ are injective;
we will identify $\pi_1(T_i,\ast)$ and $\HH_1(T_i';\Z)$ with their images in $\pi_1(\Sigma_{g,1},\ast)$
and $\HH_1(\Sigma_g;\Z)$, respectively.  Define
$K_i = [\pi,\pi] \cap \pi_1(T_i,\ast)$.
We then have the following.

\begin{lemma}
\label{lemma:birmankergen}
For $g \geq 3$, the group $[\pi,\pi]$ is generated by the $\Torelli_{g,1}$-orbits of the set $\cup_{i=1}^g K_i$.
\end{lemma}

The proof of this will have two ingredients.  The first is the following theorem of Tomaszewski.  As
notation, if $G$ is a group and $a,b \in G$, then $[a,b] := a^{-1} b^{-1} a b$ and $a^b := b^{-1} a b$.

\begin{theorem}[{Tomaszewski, \cite{Tomaszewski}}]
\label{theorem:commutatorgen}
Let $F_n$ be the free group on $\{x_1,\ldots,x_n\}$.  Then the set
$$\Set{$[x_i,x_j]^{x_i^{k_i} x_{i+1}^{k_{i+1}} \cdots x_n^{k_n}}$}{$1 \leq i < j \leq n$ and $k_m \in \Z$ for all $i \leq m \leq n$}$$
is a free basis for $[F_n,F_n]$.
\end{theorem}

\Figure{figure:birmangenset}{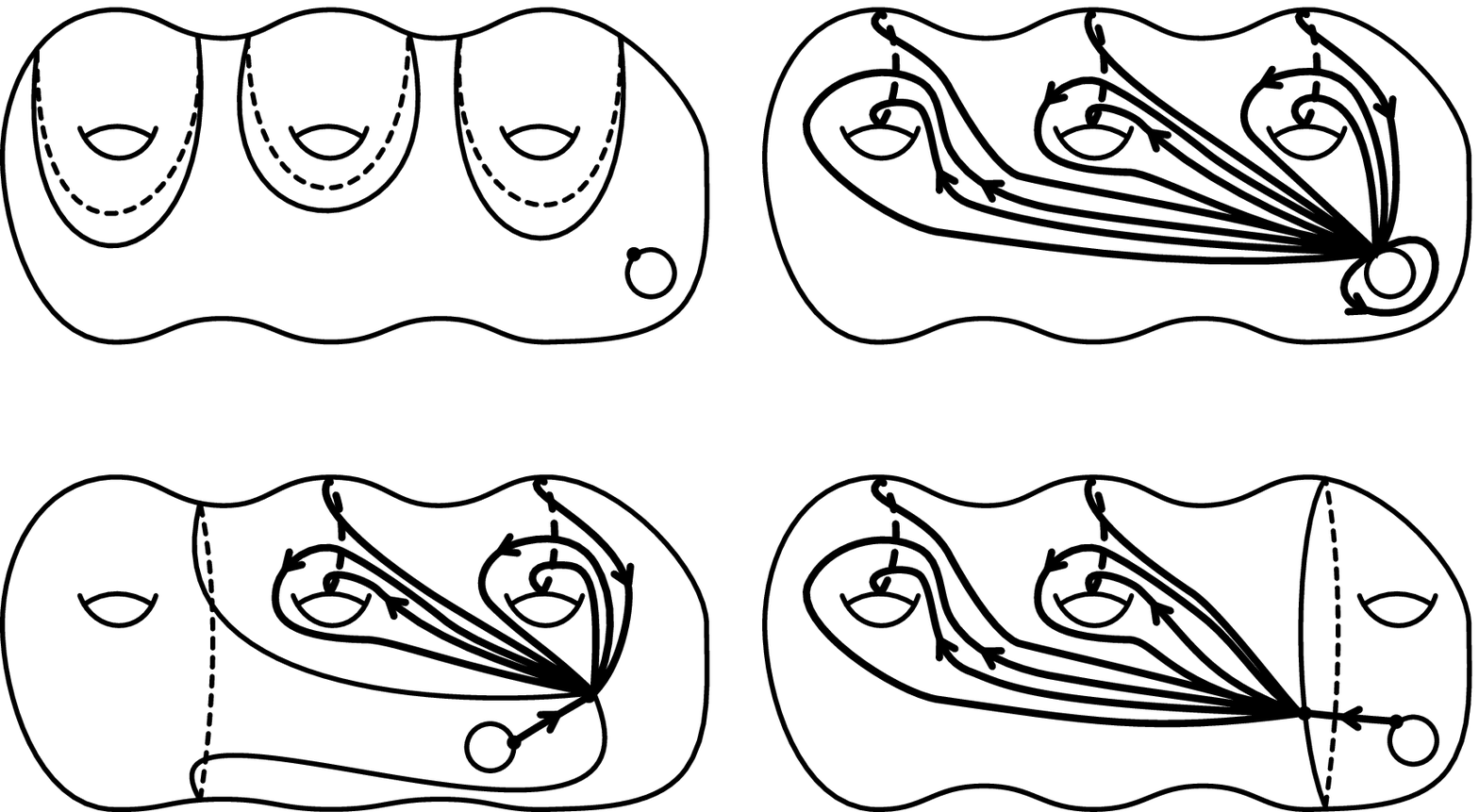}{a. The subsurfaces $T_i'$ \CaptionSpace
b. The standard basis for $\pi$ \CaptionSpace
c. The surface $X$ when $i=1$ \CaptionSpace
d. The surface $X$ when $i=g$}

\noindent
The second is the following lemma about the action of $\Torelli_{g,1}$ on $\pi$.  Choose
a standard basis $\{\alpha_1,\beta_1,\ldots,\alpha_g,\beta_g\}$ for $\pi$ (as in
Figure \ref{figure:birmangenset}.b) such that $\alpha_i$ and $\beta_i$ are freely homotopic into
$T_i'$ for $1 \leq i \leq g$.  Our proof of Lemma \ref{lemma:birmankergen} 
would be much simpler if the image of $\Mod_{g,1}$ in
$\Aut(\pi)$ contained the inner automorphisms -- since inner automorphisms act trivially on
homology, this would imply that the $\Torelli_g$-orbits of
$\Set{$[x,y]$}{$x,y \in \{\alpha_1,\beta_1,\ldots,\alpha_g,\beta_g\}$}$
generate $[\pi,\pi]$.  However, the image of $\Mod_{g,1}$ in $\Aut(\pi)$ does not
contain the inner automorphisms since $\Mod_{g,1}$
fixes the loop $\delta = [\alpha_1,\beta_1] \cdots [\alpha_g,\beta_g]$ depicted in Figure
\ref{figure:birmangenset}.b.
The following lemma is a weak replacement for this.  
\begin{lemma}
\label{lemma:torelliaction}
Let $i$ be either $1$ or $g$.  Consider $h \in \HH_1(T_i';\Z)$.  There then exists some
$w \in \Span{\alpha_i,\beta_i,\delta}$ and $f \in \Torelli_{g,1}$ such that $[w]=h$ and 
such that $f(a_j) = a_j^w$ and $f(b_j) = b_j^w$ for $1 \leq j \leq g$ with $j \neq i$.
\end{lemma}
\begin{proof}
Let $X$ be a regular neighborhood of the curves $\alpha_i \cup \beta_i \cup \partial \Sigma_{g,1}$ depicted in
Figure \ref{figure:birmangenset}.b.  Thus $X \cong \Sigma_{1,2}$, the surface $T_i'$ is homotopic into $X$, and the image of
$\pi_1(X,\ast)$ in $\pi$ is $\Span{\alpha_i,\beta_i,\delta}$.  Let 
$Y = \overline{\Sigma_{g,1} \setminus X}$, so $Y \cong \Sigma_{g-1,1}$ and $X \cap Y \cong S^1$.  
The key property of $X$ is as follows (this is where we use the assumption that $i$ is either $1$
or $g$).  There exists some $\ast' \in X \cap Y$, a properly embedded arc $\eta$ in $X$ from $\ast$ to $\ast'$,
and elements
$$\Set{$\alpha_j', \beta_j'$}{$1 \leq j \leq g$, $j \neq i$} \subset \pi_1(Y,\ast')$$
such that
$\alpha_j = \eta \cdot \alpha_j' \cdot \eta^{-1}$ and $\beta_j = \eta \cdot \beta_j' \cdot \eta^{-1}$
for $1 \leq j \leq g$ with $j \neq i$.  See Figure \ref{figure:birmangenset}.c for the case $i=1$ and Figure 
\ref{figure:birmangenset}.d for the case $i=g$.

By Lemma \ref{lemma:grafting}, there exists an oriented properly embedded arc $\eta'$ in $X$ whose
endpoints are the same as those of $\eta$ such that the homology class of $w := \eta \cdot (\eta')^{-1} \in \pi$
in $\HH_1(\Sigma_g;\Z)$ is $h$.  Observe that $w \in \Span{\alpha_i,\beta_i,\delta}$.  Also,
$$\eta' \cdot \alpha_j' \cdot (\eta')^{-1} = w^{-1} \cdot \eta \cdot \alpha_j' \cdot \eta^{-1} \cdot w = \alpha_j^w$$
for $j \neq i$, and similarly for $\beta_j$.  It is thus enough find 
some $f \in \Torelli(\Sigma_g,X)$ such that $f(\eta) = \eta'$.

The ``change of coordinates principle'' from \cite[\S 1.3]{FarbMargalitPrimer} implies that
there exists some $f' \in \Mod(\Sigma_g,X)$ such that $f'(\eta) = \eta'$.  Briefly,
an Euler characteristic calculation shows that cutting $X$ open along either $\eta$ or $\eta'$ results in
a surface homeomorphic to $\Sigma_{1,1}$.  Choosing an orientation-preserving homeomorphism
between these two cut-open surfaces and gluing the boundary components back together in an appropriate
way, we obtain some $f' \in \Mod(\Sigma_g,X)$ such that $f'(\eta) = \eta'$.  See
\cite[\S 1.3]{FarbMargalitPrimer} for more details and many other examples of arguments
of this form.  

The mapping class $f'$
need not lie in Torelli; however, it satisfies $f'([\alpha_j]) = [\alpha_j]$ and
$f'([\beta_j]) = [\beta_j]$ for $j \neq i$ and $f'(\HH_1(T_i';\Z)) = \HH_1(T_i';\Z)$.  Since the image
of $\Mod(T_i')$ in $\Aut(\HH_1(T_i';\Z)) = \Aut(\Z^2)$ is $\SL_2(\Z)$, we can choose
some $f'' \in \Mod(\Sigma_g,T_i')$ such that $f'([\alpha_i]) = f''([\alpha_i])$ and
$f'([\beta_i]) = f''([\beta_i])$.  It follows that $f := f' \cdot (f'')^{-1}$ lies in $\Torelli(\Sigma_g,X)$
and satisfies $f(\eta) = \eta'$, as desired.
\end{proof}

\begin{proof}[{Proof of Lemma \ref{lemma:birmankergen}}]
The generating set for $[F_n,F_n]$ in Theorem \ref{theorem:commutatorgen} depends on an
ordering of the generators for $F_n$.  It seems hard to prove the lemma using the
generating set corresponding to the standard ordering
$$(x_1,x_2,\ldots,x_{2g}) = (\alpha_1,\beta_1,\ldots,\alpha_g,\beta_g)$$
of the generators for $\pi \cong F_{2g}$.  However, consider the following nonstandard ordering
on the generators for $\pi$:
$$(x_1,x_2,\ldots,x_{2g}) = (\alpha_2,\beta_2,\alpha_1,\beta_1,\alpha_3,\beta_3,\alpha_4,\beta_4,\ldots,\alpha_g,\beta_g).$$
Let $S$ be the generating set for $[\pi,\pi]$ given by Theorem \ref{theorem:commutatorgen} using this
ordering of the generators.  All the elements of $S$ lie in $K_2$ except for
\begin{equation}
\label{eqn:mystery}
[\alpha_2,\zeta]^{\alpha_2^{n_2} \beta_2^{m_2} \alpha_1^{n_1} \beta_1^{m_1} \alpha_3^{n_3} \cdots \beta_g^{m_g}} \quad \text{and} \quad [\beta_2,\zeta']^{\beta_2^{m_2} \alpha_1^{n_1} \beta_1^{m_1} \alpha_3^{n_3} \cdots \beta_g^{m_g}};
\end{equation}
here $\zeta \in \{\beta_2,\alpha_1,\beta_1,\alpha_3,\ldots,\beta_g\}$ and 
$\zeta' \in \{\alpha_1,\beta_1,\alpha_3,\ldots,\beta_g\}$ and $n_i,m_i \in \Z$.  Letting $T \subset S$ be the
elements in \eqref{eqn:mystery}, we must show that every $t \in T$ can be expressed as
a product of elements in the $\Torelli_{g,1}$-orbit of the set $\cup_{i=1}^g K_i$.  
Consider $t \in T$, so either $t = [\alpha_2,\zeta]^{\alpha_2^{n_2} \beta_2^{m_2} \alpha_1^{n_1} \beta_1^{m_1} \alpha_3^{n_3} \cdots \beta_g^{m_g}}$
or $t = [\beta_2,\zeta]^{\beta_2^{m_2} \alpha_1^{n_1} \beta_1^{m_1} \alpha_3^{n_3} \cdots \beta_g^{m_g}}$.  There are two cases.

\BeginCases
\begin{case}
$\zeta \notin \{\alpha_1,\beta_1\}$.
\end{case}
We will do the case where
$t = [\alpha_2,\zeta]^{\alpha_2^{n_2} \beta_2^{m_2} \alpha_1^{n_1} \beta_1^{m_1} \alpha_3^{n_3} \cdots \beta_g^{m_g}}$; 
the other case is treated in a similar way.  Set 
$t' = [\alpha_2,\zeta]^{\alpha_2^{n_2} \beta_2^{m_2} \alpha_3^{n_3} \cdots \beta_g^{m_g}}$,
so $t' \in K_1$.  By Lemma \ref{lemma:torelliaction}, there exists some $w \in \{\alpha_1,\beta_1,\delta\}$
and $f \in \Torelli_{g,1}$ such that $[w] = [\alpha_1^{n_1} \beta_1^{m_1}]$ and such that
$f(a_j) = a_j^w$ and $f(b_j) = b_j^w$ for $j > 1$.  This implies that 
$f(t') = [\alpha_2,\zeta]^{\alpha_2^{n_2} \beta_2^{m_2} \alpha_3^{n_3} \cdots \beta_g^{m_g} w}$.  
Now, $\alpha_3^{n_3} \cdots \beta_g^{m_g} w$ and 
$\alpha_1^{n_1} \beta_1^{m_1} \alpha_3^{n_3} \cdots \beta_g^{m_g}$ are homologous, so there exists
some $\theta \in [\pi,\pi]$ such that
$\alpha_3^{n_3} \cdots \beta_g^{m_g} w \theta = \alpha_1^{n_1} \beta_1^{m_1} \alpha_3^{n_3} \cdots \beta_g^{m_g}$.
Moreover, since $w \in \Span{a_1,b_1,\delta}$ we have $\theta \in K_2$.  Observe now that
$$\theta^{-1} \cdot f(t') \cdot \theta = [\alpha_2,\zeta]^{\alpha_2^{n_2} \beta_2^{m_2} \alpha_3^{n_3} \cdots \beta_g^{m_g} w \theta}
= [\alpha_2,\zeta]^{\alpha_2^{n_2} \beta_2^{m_2} \alpha_1^{n_1} \beta_1^{m_1} \alpha_3^{n_3} \cdots \beta_g^{m_g}} = t.$$
We have thus found the desired expression for $t$.

\begin{case}
$\zeta' \in \{\alpha_1,\beta_1\}$.
\end{case}

This case is similar to Case 1.  The only difference is that the $\alpha_g^{n_g} \beta_g^{m_g}$ term
of $t$ is deleted to form $t'$ instead of the $\alpha_1^{n_1} \beta_1^{m_1}$ term.
\end{proof}

\begin{proof}[{Proof of Lemma \ref{lemma:stabgen}}]
Let $I$ be the subgroup of $\Torelli_g$ generated
by $\cup_{i=1}^g \Torelli(\Sigma_{g},S_i)$.  Using the notation of \S \ref{section:torellibdry}, there is a
surjection $\rho : \Torelli_{g,\gamma} \rightarrow (\Torelli_g)_{\gamma}$ induced by a continuous
map $\phi : \Sigma_{g,\gamma} \rightarrow \Sigma_g$.  Define $X = \phi^{-1}(S_1)$, so $X \cong \Sigma_{g-1,1}$.
Letting $\Torelli(X)$ be the Torelli group of $X$, Lemma \ref{lemma:semidirect} gives a decomposition
$\Torelli_{g,\gamma} = K_{g,\gamma} \ltimes \Torelli(X)$.
Clearly $\rho(\Torelli(X)) = \Torelli(\Sigma_g,S_1) \subset I$.  Also,
Lemma \ref{lemma:birmankergen} implies that $K_{g,\gamma}$ is generated by the $\Torelli(X)$-conjugates
of a set $S \subset K_{g,\gamma}$ such that $\rho(S) \subset I$.  We conclude that
$\rho(\Torelli_{g,\gamma}) \subset I$, as desired.
\end{proof}

\section{Proof of main theorem}
\label{section:main}

We finally prove our main theorem.  The key is the following standard lemma,
whose proof is similar to that given in \cite[(1) of Appendix to \S 3]{SerreTrees} and is thus omitted.

\begin{lemma}
\label{lemma:gensetaction}
Consider a group $G$ acting without inversions on a connected graph $X$.  Assume that
$X/G$ consists of a single edge $\overline{e}$.  Let $e$ be a lift of $\overline{e}$ to $X$ and let
$v$ and $v'$ be the endpoints of $e$.  Then $G$ is generated by $G_{v} \cup G_{v'}$.
\end{lemma}

\noindent
To apply this, we will need the following lemma.

\begin{lemma}
\label{lemma:handlegraphquotient}
Let $a,b \in \HH_1(\Sigma_g;\Z)$ satisfy $i_a(a,b)=1$.  Then $\HandleGraph / \Torelli_g$ is
isomorphic to a graph with a single edge.
\end{lemma}

\noindent
The proof is similar to the proofs of \cite[Lemma 6.2]{PutmanCutPaste} and
\cite[Lemma 6.9]{PutmanInfinite}, and is thus omitted.

\begin{proof}[{Proof of Theorem \ref{theorem:main}}]
Let $R_1',\ldots,R_g'$ and $R_{ijk}$ be the subsurfaces of $\Sigma_g$ from the introduction.
Let $\Gamma$ be the subgroup of $\Torelli_g$ generated by 
$\bigcup_{1 \leq i < j < k \leq g} \Torelli(\Sigma_g,R_{ijk})$.
Our goal is to prove that $\Gamma = \Torelli_g$.

The proof will be by induction on $g$.  The base case $g=3$ is trivial, so assume
that $g \geq 4$ and that the theorem is true for all smaller $g$ such that $g \geq 3$.  Choose simple closed
curves $\alpha$ and $\beta$ in $R_1'$ such that $i_g(\alpha,\beta)=1$.  Observe that
$R_1'$ is a closed regular neighborhood of $\alpha \cup \beta$.  Set $a = [\alpha]$ and $b = [\beta]$.  
Clearly $\Torelli_g$ acts on $\HandleGraph$ without inversions.
Lemmas \ref{lemma:handlegraphconnected} and \ref{lemma:handlegraphquotient} 
show that the action of $\Torelli_g$ on $\HandleGraph$ satisfies the other conditions of Lemma
\ref{lemma:gensetaction}.  We deduce that $\Torelli_g$ is generated by the union 
$(\Torelli_g)_{\alpha} \cup (\Torelli_g)_{\beta}$ of the stabilizer subgroups of $\alpha$ and $\beta$.

Recall that $S_i = \overline{\Sigma_g \setminus R_i'}$ for $1 \leq i \leq g$.  By Lemma \ref{lemma:stabgen},
both $(\Torelli_g)_{\alpha}$ and $(\Torelli_g)_{\beta}$ are contained in the subgroup
generated by $\cup_{i=1}^g \Torelli(\Sigma_g,S_i)$.  We must prove that $\Torelli(\Sigma_g,S_i) \subset \Gamma$ for $1 \leq i \leq g$.
We will do the case $i=g$; the other cases are similar.  We have a Birman exact sequence
$$1 \longrightarrow \pi_1(U\Sigma_{g-1}) \longrightarrow \Torelli(\Sigma_g,S_g) \longrightarrow \Torelli_{g-1} \longrightarrow 1.$$
By induction, the subset $\bigcup_{1 \leq i < j < k \leq g-1} \Torelli(\Sigma_g,R_{ijk})$ of
$\Torelli(\Sigma_g,S_g)$ projects to a generating set for $\Torelli_{g-1}$.  Also, it is clear that
the disc-pushing subgroup $\pi_1(U\Sigma_{g-1})$ of $\Torelli(\Sigma_g,S_g)$ is generated by elements that lie in
$\bigcup_{1 \leq i < j < g} \Torelli(\Sigma_g,R_{ijg})$.  We conclude that 
$\Torelli(\Sigma_g,S_g) \subset \Gamma$, as desired.
\end{proof}

\vspace{-10pt}
\begin{small}
\noindent
Andrew Putman\\
Department of Mathematics\\
Rice University, MS 136 \\
6100 Main St.\\
Houston, TX 77005\\
E-mail: {\tt andyp@rice.edu}
\end{small}

\end{document}